\documentclass[12pt]{amsart}

\usepackage[all]{xy}
\newtheorem{theorem}{Theorem}[section]
\newtheorem{lemma}[theorem]{Lemma}

\newtheorem{corollary}[theorem]{Corollary}

\theoremstyle{definition}

\theoremstyle{remark}
\newtheorem{remark}[theorem]{Remark}

\numberwithin{equation}{section}

\begin{document}

\title[Weak precompactness in Banach lattices ]
 {Weak precompactness in Banach lattices}

 \author[B. Xiang]
{Bo Xiang}
\address{School of Mathematics, Southwest Jiaotong
University, Chengdu 610031,  China}
\email{2734598188@qq.com}

\author[J.  Chen]
{Jinxi Chen }

\address{School of Mathematics, Southwest Jiaotong
University, Chengdu 610031, China}
 \email{jinxichen@swjtu.edu.cn}

\author[L. Li]
{Lei Li}
\address{Department of Mathematics, Nankai University, Tianjin 300071,  China}
\email{leilee@nankai.edu.cn}
\thanks{ The third author was partly supported  by NSFC (No.12171251).}



\begin{abstract}
	 We show that the solid hull of every weakly precompact set of a Banach lattice $E$ is weakly precompact if and only if every order interval in $E$ is weakly precompact, or equivalently, if and only if every disjoint weakly compact set is weakly precompact. Some  results on the domination property for weakly precompact positive operators are obtained. Among other things, we show that, for a pair of Banach lattices $E$ and $F$ with $E$ $\sigma$\,-Dedekind complete, every positive operator from $E$ to $F$ dominated by a weakly precompact operator is weakly precompact if and only if either the norm of $E^{\,\prime}$ is order continuous or else every order interval in $F$ is weakly precompact.
\end{abstract}

\maketitle \baselineskip 5mm

\section{Introduction}

\par Throughout this paper, $X,Y$ and $Z$ will denote real Banach spaces, $E,F$ and $G$ will denote real Banach lattices. $\it{B_X}$ denotes the closed unit ball of $X$ and $id_{X}$ denotes the identity operator on $X$. $X^{\prime}$ denotes the topological dual of $X$ and $E^{+}$ denotes the positive cone of $E$.
\par Let us recall that a bounded subset $A$ of $X$ is said to be \textit{weakly (sequentially) precompact} provided  that every sequence from $A$ has a weakly Cauchy subsequence.  Rosenthal's $\ell_{1}$-theorem asserts that a bounded set is weakly precompact if and only if it  contains no  basic sequence equivalent to the usual $\ell_{1}$ basis. It should be noted that every weakly precompact set in $X$ is relatively weakly compact if and only if $X$ is weakly sequentially complete. In terms of Dunford-Pettis operators Rosenthal characterized weakly precompact sets as follows: A subset $A$ of $X$ is weakly precompact if and only if for every Banach space $Y$ and every Dunford-Pettis operator $T:X\to{Y}$, $T(A)$ is relatively compact (\cite[p. 377]{1977Point}). From this characterization it follows immediately that the convex hull of a weakly precompact set is also weakly precompact. A bounded linear operator $T:X\to{Y}$  is called weakly precompact if $T(\it{B_X})$ is weakly precompact. Clearly, $T$ is weakly precompact if and only if $T$ does not preserve a subspace isomorphic to $\ell_{1}$. The set of all  weakly precompact operators from $X$ to $Y$ will be denoted by  $\it{WPC(X, Y)}$. It is easily verified that $\it{WPC(X,Y)}$ is a closed linear subspace of $\mathcal{L}(X,Y)$.

 \par A bounded subset $A$ in $X$ is called a \textit{limited} (resp. \textit{Dunford-Pettis}) \textit{set} if each weak$^{*}$ (resp. weakly) null sequence $(f_n)$ in $X^{\prime}$ converges uniformly to zero on $A$, i.e., $\sup_{x\in{A}}{\vert{f_n(x)}\vert}\to{0}$, or equivalently, if and only if every bounded (resp. weakly compact) operator from $X$ to $c_0$ carries $A$ to a relatively compact set \cite{2006Positive, Andrews, BD}. Clearly, every limited set is a Dunford-Pettis set. For the weak precompactness of Dunford-Pettis sets and limited sets  we have the following implications:
	\begin{center}
		limited set $\Rightarrow$ {Dunford-Pettis set } $\Rightarrow$ {weakly precompact set}
	\end{center}
   See, e.g., \cite{Andrews,1998Evaluation, BD} for details. Dually, a bounded subset $B$ of $X^{\prime}$ is called an \textit{$L$\,-set} if every weakly null sequence $(x_n)$ in $X$ tends to zero uniformly on $B$, that is, $\sup_{f\in{B}}{\vert{f(x_n)}\vert}\to{0}$ (cf. \cite{Em, 8668918320130301}). A sequence $(f_n)$ in $X^{\prime}$ is called an $L$\,-sequence if $\{f_n:n\in \mathbb{N} \}$ is an $L$\,-set.

\par Recently, in his paper \cite{Ghenciu2017A} I. Ghenciu introduced a particular class of sets that he called ``~weak limited sets". Let us recall that a bounded subset $A$ of $X$ is a \textit{weak limited set} if $T(A)$ is relatively compact for each Dunford-Pettis operator $T$ from $X$ to $c_{0}$. Among other things, he proved that  a bounded subset $A$ of $X$ is weak limited if and only if each weak$^{*}$ null $L$\,-sequence $(f_n)$ in $X^{\prime}$ converges to zero uniformly on $A$, or equivalently, if and only if it is weakly precompact (Theorem 1 and Corollary 9 of \cite{Ghenciu2017A}). Thus, a new charaterization of weak precompactness is obtained. Consequently, an operator $T:X\to{Y}$ between two Banach spaces is weakly precompact if and only if $T^{\,\prime}:{Y^{\prime}}\to{X^{\prime}}$ carries weak$^{*}$ null $L$\,-sequences to norm null sequences (see Theorem 4 and Corollary 10 of \cite{Ghenciu2017A}).
\par In this paper, one of our aims is to study weak precompactness in Banach lattices. Based upon the characterization of weak precompactness in Banach spaces due to Ghenciu \cite{Ghenciu2017A}, we consider the relationships between weak precompactness and the lattice structure of a Banach lattice.   We show that the convex solid hull of every weakly precompact set of a Banach lattice $E$ is weakly precompact if and only if every order interval in $E$ is weakly precompact, or equivalently,  if and only if every disjoint weakly compact set in $E$ is weakly precompact.  Our work improves some well known results on weak precompactness in Banach lattices.
\par Recall that in the literature the domination problem for a class  $\mathcal{C}$ of positive
operators acting between Banach lattices is stated as follows:
\begin{itemize}
  \item Let $S, T:E\rightarrow F$ be two positive operators between Banach lattices such that $0 \leq S\leq T$. Assume
that $T$ belongs to the class $\mathcal{C}$. Which conditions on $E$ and $F$ do ensure that
$S$ belongs to $\mathcal{C}$?
\end{itemize}
The literature related to the domination problems for special classes of positive operators on Banach lattices is by now very large.
\par We consider the domination problem for weakly precompact operators on Banach lattices. We prove that,  for a pair of Banach lattices $E$ and $F$ with $E$ $\sigma$\,-Dedekind complete, every positive operator from $E$ to $F$ dominated by a weakly precompact operator is weakly precompact if and only if either the norm of $E^{\,\prime}$ is order continuous or else every order interval in $F$ is weakly precompact. Furthermore, some other results on the domination property for weakly precompact  operators are obtained.

\par Our notions are standard. For the theory of Banach lattices and  operators, we refer the reader to the monographs \cite{2006Positive, 1991Banach, Wnuk}.

\section{Weak precompactness in Banach lattices}
  We start with a lemma concerning the characterization of weakly precompact sets in Banach lattices, which is due to Ghenciu \cite{Ghenciu2017A} and will be applied on several occasions in the sequel.

\begin{lemma}\cite{Ghenciu2017A}\label{Lemma 2.1} Let  $A$ be a bounded subset of a Banach space $X$. Then the following assertions are equivalent.
\begin{enumerate}
	\item $A$ is  weakly precompact.
	\item Each weak$^{*}$ null $L$\,-sequence $(f_n)$ in $X^{\prime}$ converges uniformly to zero on $A$, that is, $\mathop{sup}_{x\in{A}}{\vert{f_n(x)}\vert}\to{0}$.
\end{enumerate}
\end{lemma}
\par This lemma offers an alternative characterization of weakly precompact sets in Banach spaces. By this lemma, an operator $T:X\to{Y}$ between two Banach spaces is weakly precompact if and only if $T^{\,\prime}:{Y^{\prime}}\to{X^{\prime}}$ carries weak$^{*}$ null $L$\,-sequences to norm null sequences ( Theorem 4 and Corollary 10 of \cite{Ghenciu2017A}). In a Banach lattice $E$ every order interval is weakly precompact if and only if $\vert{f_n}\vert\xrightarrow {w^{*}} 0$ for every $L$\,-sequence $(f_n)$ of $E^{\,\prime}$ satisfing $f_{n}\xrightarrow {w^{*}} 0$. In this case, for convenience we say that the dual $E^{\,\prime}$ has \textit{$(L)$-weak$^{*}$ sequentially continuous lattice operations}.

\par Let $E$ be a Banach lattice. A bounded subset $B$ of $E^{\,\prime}$ is called  an \textit{almost $L$\,-set} if every disjoint weakly null sequence $(x_n)$ in $E$ tends to zero uniformly on $B$, that is, $\sup_{f\in{B}}{\vert{f(x_n)}\vert}\to{0}$ (see \cite{8668918320130301}). Clearly, every $L$\,-set in $E$ is an almost $L$\,-set. For a subset $S$ of a Banach lattice $E$, Sol$(S):=\{x\in E:\exists\,y\in S\,\,\textrm{with}\, \vert x\vert\leq\vert y\vert\}$ is called the solid hull of $S$. It is easily verified that the solid hull of an almost $L$\,-set is likewise an almost $L$\,-set. However, this does not necessarily hold for $L$\,-sets. We have the following lemma.
\begin{lemma}\label{Lemma 4.1} Let $E$ be a Banach lattice and let $A$ be an almost $L$\,-set in $E^{\,\prime}$.
 \begin{enumerate}
   \item If $(f_n)$ is a sequence in $Sol(A)$ satisfying ${\vert{f_n}\vert}\stackrel{w^{*}}{\longrightarrow}0$ , then $(f_n)$ is  an $L$\,-sequence.
   \item If $(f_n)$ is a disjoint sequence of Sol$(A)$, then $\vert{f_n}\vert\stackrel{w^{*}}{\longrightarrow}0$ and $(f_n)$ is  an $L$\,-sequence.
 \end{enumerate}

\end{lemma}
\begin{proof} (1) To prove $\{f_n: n\in \mathbb{N}\}$ is an $L$\,-set it is sufficient to show that ${f_n(x_n)}\to{0}$ for every weakly null sequence $(x_n)$ of $E$ (cf. Proposition 2.2 of \cite{8668918320130301}). Otherwise,  passing to a subsequence if necessary, there would exist a weakly null sequence $(x_n)$ in $E$ satisfying ${\vert{f_n}\vert}({\vert{x_n}\vert})>\vert f_{n}(x_n)\vert>{\epsilon}$ for some $\epsilon>{0}$ and for all $n$. Since ${\vert{f_n}\vert}\stackrel{w^{*}}{\longrightarrow}0$, it is easy to see that, by induction, we can find a strictly increasing subsequence $(n_{k})_{k=1}^{\infty}$ of $\mathbb{N}$ such
that $|\,f_{n_{m+1}}|(4^{m}\sum_{k=1}^{\,m}\,\vert x_{n_{k}}\vert\,)<2^{-m}$ for all $m\in\mathbb{N}$. See \cite[Theorem 2.5]{CCJ} for details. Let $$x=\sum_{k=1}^{\,\infty}2^{-k}\vert x_{n_{k}}\vert,\quad y_{m}=\left(\vert x_{n_{m+1}}\vert-4^{m}\sum_{k=1}^{\,m}\,\vert x_{n_{k}}\vert-2^{-m}x\right)^{+}.$$Then,
by Lemma 4.35 of \cite{2006Positive} $(y_{m})$ is a disjoint sequence. Now we have
\begin{eqnarray*}
    |f_{n_{m+1}}|(y_m)&=&|f_{n_{m+1}}|\left(|x_{n_{m+1}}|-4^{m}\sum_{k=1}^{\,m}\,|x_{n_{k}}|-2^{-m}x\right)^{+}\\&\geq& |f_{n_{m+1}}|\left(|x_{n_{m+1}}|-4^{m}\sum_{k=1}^{\,m}\,|x_{n_{k}}|-2^{-m}x\right)
    \\&=&|f_{n_{m+1}}|(|x_{n_{m+1}}|)-|f_{n_{m+1}}|\left(4^{m}\sum_{k=1}^{\,m}\,|x_{n_{k}}|\right)-2^{-m}|f_{n_{m+1}}|(x)\\&>&\varepsilon-2^{-m}-2^{-m}|f_{n_{m+1}}|(x).
\end{eqnarray*}
 Hence, $|f_{n_{m+1}}|(y_m)>\frac{\varepsilon}{2}$ must hold for all $m$ sufficiently large. Pick a sequence $(g_m)\subset{A}$ satisfying ${\vert{f_{n_{m}}}\vert}\le{\vert{g_m}\vert}$ for each $m$, then ${{\vert{g_{m+1}}\vert}(y_m)}>{\frac{\epsilon}{2}}$.
Thus there exists a sequence $(u_m)$ in $E$ such that  ${\vert{u_m}\vert}\le{y_m}$ and ${\vert{g_{m+1}(u_m)}\vert}>{\frac{\epsilon}{2}}$ must hold for  $m$ sufficiently large. Clearly, $(u_m)$ is also disjoint. Since ${\vert{u_m}\vert}\le{y_m}\le{\vert{x_{n_{m+1}}}\vert}$ and $x_{n_{m+1}} \stackrel{w}{\longrightarrow}0$, from Theorem 4.34 of \cite{2006Positive} it follows that $u_m\stackrel{w}{\longrightarrow}0$. Since $A$ is an almost $L$\,-set of $E^{\,\prime}$, we have ${\vert{g_{m+1}(u_m)}\vert}\le{\sup_{f\in{A}}{\vert{f}(u_m)\vert}\to{0}}$ , which is impossible. This proves that $\{f_n: n\in \mathbb{N}\}$ is an $L$\,-sequence of $E^{\,\prime}$.
\par (2) Let $(f_n)$ be an arbitrary disjoint sequence from Sol$(A)$.  Then  $\vert{f_n}\vert\stackrel{w^{*}}{\longrightarrow}0$ in $E^{\,\prime}$ (cf.  \cite[ Corollary 3.3]{FaOu2021}). From (1) it follows that $(f_n)$ is an $L$\,-sequence.
\end{proof}

\par  The convex hulls of  weakly precompact sets in  Banach spaces  are all  weakly precompact. However, a weakly precompact set in a Banach lattice does not necessarily have a weakly precompact solid hull. By Proposition 2.5.12 ii) of \cite{1991Banach} we know that every weakly precompact subset of a Banach lattice $E$ has a weakly precompact solid hull if either the norm on $E$ is order continuous or else $E^{\,\prime}$ has a weak order unit. Wickstead showed that every weakly relatively compact subset of $E^{+}$ has a weakly relatively compact solid hull if and only if $E$ has order continuous norm, that is, every order interval in $E$ is weakly compact (\cite{Wickstead}, cf. \cite[Theorem 4.39]{2006Positive}). For the solid hulls of weakly precompact sets we have the following result.
\begin{theorem}\label{solid hull}
Let $E$ be a Banach lattice. Then the solid hull of every weakly precompact set of $E$ is weakly precompact if and only if every order interval in $E$ is weakly precompact, that is, $E^{\,\prime}$ has $(L)$-weak$^{*}$ sequentially continuous lattice operations.
\end{theorem}
\begin{proof}
If the solid hull of  every weakly precompact set of $E$ is weakly precompact, then, for each $x\in E^{+}$, $[-x,x]=\textmd{Sol}\{x\}$ is obviously weakly precompact.
\par For the converse, let $A$ be a weakly precompact set of $E$. To prove that Sol$(A)$ is also weakly precompact we have to show that $$\sup\{|f_{n}(z)|:z\in \textmd{Sol}(A)\}\rightarrow0$$for every weak$^{*}$ null $L$\,-sequence $(f_n)$ of $E^{\,\prime}$. Otherwise, there would exist a weak$^{*}$ null $L$\,-sequence $(f_n)$ of $E^{\,\prime}$ and a sequence $z_{n}\in \textmd{Sol}(A)$ such that $$\epsilon<|f_{n}(z_{n})|\le|f_{n}|(|z_{n}|)$$ for some $\epsilon>0$ and for all $n\in \mathbb{N}$. Find $x_{n}\in A$ satisfying $|z_n|\le|x_{n}|$ for each $n$. Then we have $|f_{n}|(|x_{n}|)>\epsilon$ for each $n\in \mathbb{N}$. By the Riesz-Kantorovich formula, for each $n$, there exists $g_n\in E^{\,\prime}$ such that $|g_n|\le|f_n|$ and $|g_{n}(x_{n})|>\epsilon$. Since $(f_n)$ is a weak$^{*}$ null $L$\,-sequence in $E^{\,\prime}$, by hypothesis, we have $|f_n|\xrightarrow{w^*}0$, hence $|g_n|\xrightarrow{w^*}0$. By Lemma \ref{Lemma 4.1}(1) we know that $(g_n)$ is also a weak$^{*} $ null $L$\,-sequence in $E^{\,\prime}$. Therefore, from the weak precompactness of $A$ it follows that $$\epsilon<|g_{n}(x_{n})|\le \sup\{|g_n(x)|:x\in A\}\to 0\,\,\,  (n\to \infty)$$This is impossible, and the proof is finished.
\end{proof}
\par In his monograph W. Wnuk \cite{Wnuk} called a bounded subset $A$  of a Banach lattice $E$ a disjoint weakly compact set whenever every disjoint sequence $(x_n)\subset \textmd{Sol}(A)$ is weakly null. Note that every weakly precompact set in a Banach lattice is disjoint weakly compact \cite[Proposition 2.5.12 iii) ]{1991Banach}. For a non-empty bounded subset $A$ of $E$, we define a lattice seminorm
$$\rho_{A}(f)=\sup\{\langle|f|,\,|x|\rangle:x\in A\},\qquad f\in E^{\,\prime}$$
By the Riesz-Kantorovich formula it is easy to verify that $$\rho_{A}(f)=\sup\{|\langle f,\,z\rangle|:z\in \textmd{Sol}(A)\}.$$If $A$ is a weakly precompact set in $E$, then $\rho_{A}(g_n)\to{0}$ for every order bounded disjoint sequence $(g_n)\subset E^{\,\prime}$ \cite[Corollary 2.5.2]{1991Banach}. On the other hand, the reverse implication also holds if either the norm on $E$ is order continuous or  $E^{\,\prime}$ has a weak order unit. \cite[Corollary 2.5.9]{1991Banach}. This is improved by the following result.
\begin{theorem}\label{supplement 1} Let $E$ be a Banach lattice. Then the following  assertions are equivalent.
\begin{enumerate}

  \item For a non-empty norm bounded set $A\subset E$, $A$ is weakly precompact if and only if $\rho_{A}(g_n)\to{0}$ holds for every order bounded disjoint sequence $(g_n)\subset E^{\,\prime}$.
  \item For a non-empty norm bounded set $A\subset E$, $A$ is weakly precompact if and only if $A$ is disjoint weakly compact.
 \item Every order interval in $E$ is weakly precompact.
\end{enumerate}
\end{theorem}
\begin{proof}
(1)$\Rightarrow$(3) Let $x\in E^{+}$. For every order bounded disjoint sequence $(g_n)\subset E^{\,\prime}$, we have $|g_n|\xrightarrow{w}0$. Therefore, $\rho_{\textmd{Sol}\{x\}}(g_n)=|g_n|(x)\to{0}$. By our hypothesis, $[-x,\,x]=\textmd{Sol}\{x\}$ is weakly precompact.
\par (1)$\Leftrightarrow$(2) It should be first noted that for a non empty subset of a Banach lattice $E$ the following two conditions are equivalent:
\begin{itemize}
  \item Each order bounded disjoint sequence in $E^{\,\prime}$ converges uniformly to zero on $\textmd{Sol}(A)$.
  \item Each disjoint sequence in $\textmd{Sol}(A)$ (absolutely) weakly converges to zero.
\end{itemize}
 The equivalence follows easily from a general result due to O. Burkinshaw and P. G. Dodds (cf. \cite[Theorem 5.63]{2006Positive}). Then this observation implies that the proof of (1)$\Leftrightarrow$(2) is easy.
\par (3)$\Rightarrow$(2) Assume that every order interval in $E$ is weakly precompact. We only need to prove that the set $A$ is weakly precompact if $A$ is disjoint weakly compact. We can assume without loss of generality that $A$ is solid.  So, let $A$ be a solid subset of $E$ such that each disjoint sequence in $A$ (absolutely) weakly converges to zero.
 \par Assume by way of contradiction  that $A$ is not weakly precompact. Then, there would exist a weak$^{*}$ null $L$\,-sequence $(f_n)$ in $E^{\,\prime}$ and a sequence $(x_n)$ in $A$ satisfying $${\vert{f_n}\vert}({\vert{x_n}\vert})>\vert f_{n}(x_n)\vert>{\epsilon}$$ for some $\epsilon>{0}$ and for all $n$. Note that $|f_n|\xrightarrow{w^*}0$ since every order interval in $E$ is weakly precompact. Also, from Lemma \ref{Lemma 4.1} it follows that $(|f_n|)_{1}^{\infty}$ is itself an $L$\,-sequence.   By applying the construction technique of disjoint sequences in the proof of Lemma \ref{Lemma 4.1}, we can find a subsequence $(f_{n_{m}})$ of $(f_{n})$ and  a disjoint sequence $(y_m)_{1}^{\infty}\subset A^{+}$ such that  $|f_{n_{m+1}}|(y_m)>\frac{\varepsilon}{2}$ must hold for all $m$ sufficiently large. By our hypothesis, we have $y_m\xrightarrow{w}0$. However, by definition, $$\frac{\varepsilon}{2}<|f_{n_{m+1}}|(y_m)\to{0} \qquad(n\to{\infty})$$ since $(|f_{n_{m+1}}|)_{m=1}^{\infty}$ is an $L$\,-sequence. This is impossible. The set $A$ must be a weakly precompact set, as desired.
\end{proof}
\par Now we will apply the previous theorem to deduce a Kadec-Pelczynski type result, which generalizes Proposition 2.5.13 of \cite{1991Banach}. Here, the proof is a little more direct.
\begin{corollary}\label{supplement 2} Let $E$ be a Banach lattice for which every order interval is weakly precompact. Then for a bounded solid subset $A$ of $E$ the following two assertions are equivalent.
\begin{enumerate}
  \item $A$ is weakly precompact.
  \item $A$ contains no disjoint sequence equivalent to the unit vector basis of $\ell_{1}$.
\end{enumerate}
\end{corollary}
\begin{proof}
$(1)\Rightarrow(2)$  follows immediately from Rosenthal's $\ell_{1}$-theorem.
\par $(2)\Rightarrow(1)$ Assume by way of contradiction that the set $A$ fails to be weakly precompact. By Theorem \ref{supplement 1} there would exist a disjoint sequence $(x_{n})_{1}^{\infty}\subset A^{+}$ which fails to be weakly null. Then there exists some $f\in (E^{\,\prime})^{+}$, passing to a subsequence if necessary, such that $f(x_{n})\geq1$ for all $n\in\mathbb{N}$. From this it follows that
 $$\sum_{i=1}^{n}|\alpha_{i}|\leq \left\langle f, \sum_{i=1}^{n}|\alpha_{i}|x_{n}\right\rangle=\left\langle f, \left|\sum_{i=1}^{n}\alpha_{i}x_{n}\right|\right\rangle\leq\|f\|\left\|\sum_{i=1}^{n}\alpha_{i}x_{n}\right\|$$holds for every $n\in\mathbb{N}$ and every choice of scalars $\alpha_{1}, \cdot\cdot\cdot, \alpha_{n}$. This implies that $(x_{n})_{n}^{\infty}$ is equivalent to the unit vector basis $(e_{n})_{1}^{\infty}$ of $\ell_{1}$.
\end{proof}

\section{ domination by weakly precompact operators}

 \par A continuous operator $T:E\to{X}$ from a Banach lattice $E$ into a Banach space $X$ is said to be \textit{order weakly precompact} if $T[-x,x]$ is weakly precompact in $X$ for all $x\in{E^{+}}$. The set of all order weakly precompact operators from $E$ to $X$ will be denoted by $\it{oWPC(E,X)}$. It is easily verified that $\it{oWPC(E,X)}$ is a norm closed vector subspace of $\it{\mathcal{L}(E,X)}$. Recall that an operator $T:E\to{X}$ is\textit{ AM-compact} (resp. \textit{order weakly compact}) if $T[-x,x]$ is relatively (resp.  weakly) compact  in $X$ for all $x\in{E^{+}}$.

\begin{theorem}\label{Theorem 3.3} Let $T:E\to{X}$ be a bounded linear operator from a Banach lattice $E$ into a Banach space $X$. Then the following assertions are equivalent.
\begin{enumerate}
	\item $T$ is an order weakly precompact operator.
	\item For each Dunford-Pettis operator $S$ from $X$ into arbitrary Banach space $Y$, $ST$ is AM-compact.
	\item For each Dunford-Pettis operator $S$ from $X$ into $c_{0}$, $ST$ is AM-compact.
	\item For each weak$^{*}$ null $L$-sequence $(f_n)$ of $X^{\prime}$ we have ${\vert{T^{\,\prime}f_n}\vert} \stackrel{w^{*}}{\longrightarrow}0$ in $E^{\prime}$.
\end{enumerate}
\end{theorem}
\begin{proof} (1)$\Rightarrow$(2)$\Rightarrow$(3) Obvious.
\par (3)$\Rightarrow$(4) Let $(f_n)$ be a weak$^{*}$ null $L$\,-sequence of $X^{\prime}$, we consider the operator $S:X\to{c_{0}}$ defined by $S(x)=(f_n(x))_{n=1}^{\infty}$ for all $x\in{X}$. Let $(x_n)$ be an arbitrary weakly null sequence in $X$. Since $(f_n)$ be an $L$\,-sequence of $X^{\prime}$, then
$$\begin{aligned}
	\lim_{n}\|S(x_n)\|=\lim_{n}\sup_{i}{\vert{f_{i}(x_n)}\vert}=0
\end{aligned}$$
Thus $S$ is a Dunford-Pettis operator. By assumption $ST([-x, x])$ is relatively compact subset of $c_{0}$ for all $x\in{E^{+}}$. It follows that
$$\begin{aligned}
	\vert{T^{\prime}f_n}\vert(x)=\sup{\{\vert{T^{\,\prime}f_n}(y)\vert}:y\in{[-x,x]}\}=\sup{\{\vert{f_n}(z)\vert}:z\in{T[-x,x]}\}\to{0}.
\end{aligned}$$
\par (4)$\Rightarrow$(1) For each $x\in{E^{+}}$, we have to show that $T[-x,x]$ is a weakly precompact set of $X$. To this end, by Lemma \ref{Lemma 2.1} let $(f_n)$ be an arbitrary weak$^{*}$ null $L$-sequence of $X^{\prime}$, then
$\sup{\{\vert{f_n}(y)\vert}:y\in{T[-x,x]}\}={{\vert{T^{\,\prime}f_n}\vert}}(x)\to{0}.$ It follows that $T[-x,x]$ is a weakly precompact set of $X$ for each $x\in{E^{+}}$, and hence $T$ is an order weakly precompact operator.
\end{proof}
\par Let us recall that a  bounded subset $A$ of a Banach lattice $E$ is called  \textit{$L\,$-weakly compact} if $\|x_n\|\rightarrow0$ for every disjoint sequence $(x_{n})$ contained in the solid hull of $A$ (cf. \cite[Definition 3.6.1]{1991Banach}). Every $L$\,-weakly compact set is relatively weakly compact set. Every relatively weakly compact subset of $E$ is $L$\,-weakly compact if and only if $E$ has the positive Schur property, that is,  every weakly null sequence with positive terms is norm null (cf. \cite[Corollary 3.6.8]{1991Banach}).

\begin{theorem}\label{Theorem 3.5}
	For a Banach lattice $E$ the following statements are equivalent.
\begin{enumerate}
  \item The norm of $E^{\,\prime}$  is order continuous.
  \item Every almost $L$\,-set in $E^{\,\prime}$ is $L$\,-weakly compact.
  \item Every order weakly precompact operator from $E$ to an arbitrary Banach space $X$ is a weakly precompact operator.
  \item Every order weakly compact operator from $E$ to an arbitrary Banach space $X$ is a weakly precompact operator.
\end{enumerate}
\end{theorem}
\begin{proof}
(1)$\Leftrightarrow$(2) is  Lemma 2.9 of \cite{CCJ}.
\par(2)$\Rightarrow$(3) Let $T:E\to{X}$ be an order weakly precompact operator.  To prove that $T$ is weakly precompact we have to prove that $\|T^{\,\prime}f_n\|\to{0}$  for every  $L$\,-sequence $(f_n)$ in $X^{\prime}$ satisfying $f_n\stackrel{w^{*}}{\longrightarrow}0$. To this end,  let $(f_{n})\subset{X^{\prime}}$ be an arbitrary weak$^{*}$ null $L$\,-sequence. By a result due to Dodds and Fremlin  it suffices to prove that $\vert{T^{\,\prime}f_n}\vert\stackrel{w^{*}}\longrightarrow{0}$  and $T^{\,\prime}f_n(x_n)\to{0}$ for every norm bounded disjoint sequence $(x_n)$ in $E^{+}$ (\cite{1979Compact}, cf. \cite[Proposition 2.3.4]{1991Banach}).
Since $T$ is order weakly precompact, by Theorem \ref{Theorem 3.3} we have ${\vert{T^{\,\prime}f_n}\vert}\stackrel{w*}{\longrightarrow}0$.
 It is clear that $A=\{T^{\,\prime}f_{n}:{n}\in{\mathbb{N}}\}$ is an $L$\,-set, hence an almost $L$\,-set of $E^{\,\prime}$. Then, by our hypothesis, $A$ is  $L$\,-weakly compact. Let $\rho_{A}(x)=\sup\{\vert T^{\,\prime} f_{n}\vert (\vert{x}\vert): n\in\mathbb{N}\}$. Let $(x_{n})$ be an arbitrary norm bounded disjoint sequence of $E^{+}$. It follows from Proposition 3.6.3 of \cite{1991Banach} that
	\begin{center}
	$\vert{T^{\,\prime}f_{n}(x_{n})}\vert\leq{\vert{T^{\,\prime}f_{n}}\vert}{(\vert x_{n}\vert)}\leq\rho_{A}(x_{n})\to{0}$.
	\end{center}

Now, we have $\|T^{\,\prime}f_{n}\|\to{0}$, and hence $T$ is weakly precompact.
\par (3)$\Rightarrow$(4) Obvious.
\par(4)$\Leftrightarrow$(1) is contained in \cite[Theorem 3.4.18]{1991Banach}.
\end{proof}

\par We are now in  a position to give some sufficient and necessary conditions for which a positive operator dominated by a weakly precompact operator is always weakly precompact.

\begin{theorem}\label{theorem 4.0} Let $E$ and $F$ be two Banach lattices such that every order inteval in $F$ is weakly precompact (in particular, either $F$ has order continuous norm or  $F^{\,\prime}$ has a weak order unit). Then for all operators $S, T:E\to{F}$ such that $0\le{S}\le{T}$ and $T$ is weakly precompact, the operator $S$ is weakly precompact.
\end{theorem}
\begin{proof}By Theorem \ref{solid hull}, the solid hull of every weakly precompact set in $F$ is weakly precompact if every order inteval in $F$ is also weakly precompact. Note that $S(B_E)$ is contained in $\textmd{Sol}(T(B_E))$. Since $T$ is weakly precompact, it follows that $S(B_E)$ is weakly precompact.
\end{proof}

\par If the unit ball $B_E$ of a Banach lattice $E$ is weakly precompact, then $E$ does not contain an isomorphic copy of $\ell_1$. This implies that $E^{\,\prime}$ is a KB-space, that is, the norm on $E^{\,\prime}$ is order continuous (cf. \cite[Theorem 4.69]{2006Positive}). On the other hand, if $E^{\,\prime}$ is a KB-space with a weak order unit, then from Theorem 2.5.9 of \cite{1991Banach} it follows that $B_E$ is weakly precompact. Of course, in this case every bounded linear operator from or into $E$ is weakly precompact. However, there exists a $\sigma$\,-Dedekind complete Banach lattice (e.g., $\ell_\infty$ ) whose dual has an order continuous norm but does not have a weak order unit.

\begin{theorem}\label{theorem 4.4}Let $E$ and $F$ be two Banach lattices such that  the norm of $E^{\,\prime}$ is order continuous. If, in addition, either $E$ is $\sigma$\,-Dedekind complete or $E^{\,\prime}$ has a weak order unit, then for all operators $S,T:E\to{F}$ such that $0\le{S}\le{T}$ and $T$ is weakly precompact, the operator $S$ is weakly precompact.
\end{theorem}
\begin{proof}
 We only need to consider the nontrivial case when $E$ is  $\sigma$\,-Dedekind complete and $E^{\,\prime}$ has order continuous norm.  First we recall that every weakly precompact operator from a $\sigma$\,-Dedekind complete Banach lattice $E$ to a Banach lattice $F$ is order weakly compact (cf. \cite[Proposition 3.4.15]{1991Banach}). Since $T$ is weakly precompact, $T$ is order weakly compact. Hence $S$, dominated by $T$, is also order weakly compact. Since the norm on $E^{\,\prime}$ is order continuous, it follows from Theorem \ref{Theorem 3.5} that $S$ is weakly precompact.
\end{proof}
\par As an immediate consequence we give a result on the domination problem for order weakly precompact operators.

\begin{theorem}\label{Theorem 4.11}
 Let $E$ and $F$ be two Banach lattices such that $E$ is $\sigma$\,-Dedekind complete. If a positive operator $S:E\to{F}$ is dominated by an order weakly precompact operator, then $S$ is order weakly precompact.
 \end{theorem}
 \begin{proof}Let $x\in{E^{+}}$. Then the principal ideal $E_{x}$ generated by $x$ in $E$  is an \textit{AM}-space with respect the norm ${{\|{y}\|}_{\infty}}=\inf{\{{\lambda>0}: {{\vert{y}\vert}\le{\lambda{x}}}\}}$ $(y\in E_{x})$. Also $E_{x}$ has $x$ as order unit and $[-x,x]$ as its closed unit ball.   Note that $E_{x}$ is  $\sigma$\,-Dedekind complete since $E$ is also $\sigma$\,-Dedekind complete. It is easy to see that $S:E\to{F}$ is order weakly precompact if and only if for each $x\in{E^{+}}$ the restriction operator $S:(E_{x},\|\cdot\|_{\infty})\to{F}$ is a  weakly precompact operator.
 \par Let $T$ be an order weakly precompact operator such that $0\le{S}\le{T}:E\to{F}$, and  let $x\in{E^{+}}$. Clearly, $0\le S\le T:(E_{x},\|\cdot\|_{\infty})\to{F}$. Since $T:{(E_{x},{\|{\cdot}\|}_{\infty})}\to{F}$ is weakly precompact,  and the norm of $(E_{x})^{\prime}$ is order continuous, it follows from Theorem \ref{theorem 4.4} that $S:{(E_{x},{\|{\cdot}\|}_{\infty})}\to{F}$ is likewise weakly precompact, that is, $S$ is order weakly precompact, as desired.
\end{proof}

\par The following result states that in order that every positive operator $S:E\to F$  dominated by a weakly precompact operator is weakly precompact, it is necessary that either the norm of $E^{\,\prime}$ be  order continuous or every order interval in $F$ be weakly precompact.

\begin{theorem}\label{Theorem 4.5}
Let $E$ and $F$ be two Banach lattices. If for all operators $S,T:E\to{F}$ such that $0\le{S}\le{T}$ and $T$ is weakly precompact, the operator $S$ is weakly precompact, then either the norm of $E^{\,\prime}$ is  order continuous or else every order interval in $F$ is weakly precompact.
\end{theorem}
\begin{proof}Assume by way of contradiction that neither is the norm of $E^{\,\prime}$  order continuous nor are the lattice operations in $F^{\,\prime}$  $(L)$\,-weak$^{*}$ sequentially continuous. We have to construct a positive operator which is not weakly precompact while it is dominated by a weakly precompact operator. Our construction technique is inspired by \cite{Wickstead, Wickstead 2}
 \par If the norm of $E^{\,\prime}$ is not order continuous,  there would exist an order bounded disjoint sequence $(\phi_n)\subset[0, \psi]$ in $(E^{\,\prime})^{+}$ satisfying $\|\phi_n\|=1$ for all $n$. We define the operator $S_{1}:E\to{\ell_{1}}$ as
	\begin{center}
		$S_{1}x=(\phi_n(x))_{n=1}^{\infty}$
	\end{center}
for all $x\in{E}$ since ${\sum_{n=1}^{\infty}{\vert{\phi_n(x)}\vert}}\le{\sum_{n=1}^{\infty}{\phi_n{(\vert{x}\vert)}}}={\psi(\vert{x}\vert)}$ holds for each $x\in{E}$. Clearly $S_{1}$ is positive. If the lattice operations in $F^{\,\prime}$ are not $(L)$\,-weak$^{*}$ sequentially continuous, there would exist a weak$^{*}$ null $L$\,-sequence $(f_n)\subset{F^{\,\prime}}$ such that $(\vert{f_n}\vert)$ is not weak$^{*}$ null. Therefore, by passing to a subsequence if necessary, we can find some $\epsilon>{0}$ and some $y\in{F^{+}}$ satisfying
	${\vert{f_n}\vert}(y)>{2\epsilon}$ for all ${n}\in{\mathbb{N}}$. By the Riesz-Kantorovich formula we can find  $y_n\in [0,y]$ satisfying  ${\vert{{f_n}(y_n)}\vert}>{\epsilon}$ for each $n\in \mathbb{N}$. We let $S_{2}:{\ell_{1}}\to{F}$ be defined by
		$${S_{2}((\lambda_{n}))}={\sum_{n=1}^{\infty}{\lambda}_{n}{y_n}}$$ for all ${({\lambda}_{n}})\in{\ell_{1}}$.
\par Now we define two positive operators $S,T:E\to{F}$ by
	\begin{center}
		$S(x)=S_{2}S_{1}(x)={{\sum_{n=1}^{\infty}}{\phi_n(x)}{y_n}}$\,\, and \,\,\,$T(x)=\psi(x)y$.
	\end{center}
for all $x\in E$. Obviously $T$ is compact, and hence $T$ is weakly precompact. Also, it is easy to verify that $0\le{S}\le{T}$. For each $n\in \mathbb{N}$, let $S^{\prime}_{n}:F^{\prime}\to{E^{\prime}}$ be defined as
      $$S^{\prime}_{n}f=\sum_{i=1}^{n}f(y_{i})\phi_i=\langle f, \sum_{i=1}^{n}y_{i}\otimes \phi_i\rangle$$for all $f\in F^{\prime}$. Clearly $0\le S^{\prime}_{n}\uparrow_{n}$ and $S^{\prime}_{n}f(x)\uparrow_{n} S^{\prime}f(x)$ for all $0\le f\in F^{\prime}$ and all $x\in E^{+}$. It follows from Theorem 1.19 of \cite{2006Positive} that $ 0\le S^{\prime}_{n}\uparrow_{n} S^{\prime}$ in $\mathcal{L}_{b}(F^{\prime}, E^{\prime})$, the Dedekind complete Riesz space of all order bounded operator from $F^{\,\prime}$ to $E^{\,\prime}$. Then $S^{\prime}_{n}f\xrightarrow {o}S^{\prime}f$, hence $\vert S^{\prime}_{n}f\vert\xrightarrow {o}\vert S^{\prime}f\vert\,\, (n\rightarrow\infty)$ for all $f\in F^{\prime}$. Since $(\phi_n)\subset[0, \psi]$ in $(E^{\,\prime})^{+}$ is  disjoint, for all $f\in F^{\prime}$ we have
\begin{eqnarray*}
      \vert S^{\prime}f\vert= o-\lim_{n\rightarrow\infty}\vert S^{\prime}_{n}f\vert&=&o-\lim_{n\rightarrow\infty}\sum_{i=1}^{n}\vert f(y_{i})\vert \phi_{i}\\&=&o-{{\sum_{i=1}^{\infty}}{\vert{f(y_{i})}\vert}{\phi_{i}}}\ge{{\vert{f(y_n)}\vert}{\phi_n}}
\end{eqnarray*}for all $n\in \mathbb{N}$. (Here, the symbols `\,$\xrightarrow {o}$\,' and `\,$o-\lim$\,' stand for order convergence and order limit, respectively.) In particular,
${\vert{S^{\prime}(f_n)}\vert}\ge{{\vert{f_n(y_n)}\vert}\,{\phi_n}}$ for all $n$. Therefore,
	\begin{center}
		$\|S^{\prime}f_n\|\ge{{\vert{f_n(y_n)}\vert}{\|\phi_n}\|}>{\epsilon}$
	\end{center}
	for all $n$. This implies that $S$ is not weakly precompact. The proof is finished.
 \end{proof}

\begin{remark}\label{Remark 4.6} Let $E$ and $F$ be a pair of Banach lattices. A well known result due to A. W. Wickstead \cite{Wickstead} tells us that every positive operator from  $E$ to $F$ dominated by a weakly compact positive operator is weakly compact if and only if either $E^{\,\prime}$ or $F$ has an order continuous norm. Note that the order continuity of the norm on $F$ is equivalent to the weak compactness of order intervals in $F$.
Naturally, for the domination problem of weakly precompact positive operators the authors pose the following question:

\noindent\textbf{Question.} 	\textit{Are the following assertions are equivalent for a pair of Banach lattices $E$ and $F$?}
	\begin{enumerate}
		\item \textit{Every positive operator $S:E\to F$  dominated by a weakly precompact operator is weakly precompact.}
		\item \textit{One of the following two conditions holds:}
		\begin{enumerate}
			\item \textit{The norm of $E^{\,\prime}$ is order continuous.}
			\item \textit{Every order interval in $F$ is weakly precompact.	}
		\end{enumerate}	
	\end{enumerate}
Clearly, the answer to this question is affirmative when $E$ is $\sigma$\,-Dedekind complete. For the general case we are not able to answer it and have to leave it open.
\end{remark}

\begin{theorem}\label{Theorem 4.7}
  Let $E,\,F$ be Banach lattices and $X$ be a Banach space. Consider the scheme of operators
$E\stackrel{S_{1}}{\longrightarrow}F\stackrel{S_{2}}{\longrightarrow}X$. If the positive operator $S_{1}$ is dominated by a weakly precompact operator and the operator $S_{2}$ is order weakly compact, then $S_{2}\circ{S_{1}}$ is weak precompact.
\end{theorem}
\begin{proof} Since $S_{2}$ is order weakly compact, we know $S_2$ admits a factorization through  a Banach lattice $H$ with an order continuous norm

$$\xymatrix
 {E\ar[r]^{S_{1}}&F\ar[rr]^{S_{2}}\ar[dr]_{Q}&&X \\
&&H\ar[ur]_{T}&
}$$
such that $Q:F\to{H}$ is a lattice homomorphism (cf.\,\cite[Theorem 5.58]{2006Positive}\,). Since the positive operator $QS_{1}:E\to{H}$ is still dominated by a weakly precompact operator and $H$ has an order continuous norm,  it follows from Theorem \ref{theorem 4.0} that $QS_{1}$ is weakly precompact. Hence $S_{2}S_{1}=TQS_{1}$ is weakly precompact.
\end{proof}
\par
Note, once more, that every weakly precompact operator  from a $\sigma$\,-Dedekind complete Banach lattice  to another Banach lattice is order weakly compact \cite[Proposition 3.4.15]{1991Banach}. As immediate consequences of the preceding theorem we have the following results.
\begin{corollary}\label{Corollary 4.8}
 Let $E, F$ and $G$ be Banach lattices such that $F$ is $\sigma$\,-Dedekind complete. Let $S_{1},T_{1}:E\to{F}$ and $S_{2},T_{2}:F\to{G}$ be operators such that $0\le{S_{i}}\le{T_{i}}$ and each $T_{i}$ is weakly precompact , $i=1,2$. Then $S_{2}S_{1}$ is weakly precompact operator.
\end{corollary}
\begin{corollary}\label{Corollary 4.9}
 If a positive operator $S$ on a $\sigma$\,-Dedekind complete Banach lattice is dominated by a weakly precompact positive operator, then $S^{2}$ is also a weakly precompact operator.
\end{corollary}
\par
For positive orthomorphisms dominated weakly precompact operators we have the analogue of Theorem 5.21 and Theorem 5.34 of \cite{2006Positive} as follows, and the proof is very similar.
\begin{theorem}\label{Theorem 4.10}
 A positive orthomorphism on a Dedekind complete Banach lattice dominated by a weakly precompact operator is itself weakly precompact.
 \end{theorem}
\begin{proof} Let $E$ be a Dedekind complete Banach lattice and $S,T:E\to{E}$ be two positive operators such that ${0}\le{S}\le{T}$ holds. Assume that $S$ is an orthomorphism and that $T$ is weakly precompact. We can assume without loss of generality that ${0}\le{S}\le{id_{E}}$ holds. Now let $\epsilon>{0}$. Then by Theorem 2.9 of \cite{2006Positive}, there exists an ${id_{E}}$\,-step function $\sum_{i=1}^{n}{\alpha_{i}P_{i}}$ with ${0}\le{\sum_{i=1}^{n}{\alpha_{i}P_{i}}}\le{S}$ and ${\|{S-{\sum_{i=1}^{n}{\alpha_{i}P_{i}}}}\|}<{\epsilon}$. Here  $P_{i}$ is an order projection on $E$ and ${\alpha_{i}}>{0}$ holds for each $i$. Since ${0}\le{P_{i}}\le{\frac{1}{\alpha_{i}}{T}}$, from  Corollary \ref{Corollary 4.9} it follows that each $P_{i}={P_{i}^{2}}$ is weakly precompact, and so $\sum_{i=1}^{n}{\alpha_{i}P_{i}}$ is also weakly precompact. It is easy to see that $S$ must be a weakly precompact operator since $WPC(E, E)$, the collection of all weakly precompact operators on $E$,  is a closed linear subspace of $\mathcal{L}(E,E)$.
\end{proof}



\end{document}